\documentclass[11pt]{amsart}
\usepackage{graphicx}
\usepackage{stmaryrd}
\usepackage{amssymb}
\usepackage{amsthm}
\usepackage{dsfont}
\usepackage{hyperref}

\makeatletter
\@namedef{subjclassname@2010}{%
  \textup{2010} Mathematics Subject Classification}
\makeatother

\newtheorem{cor}{Corollary}
\newtheorem{lem}{Lemma}
\newtheorem{rem}{Remark}
\newtheorem{prop}{Proposition} 

\numberwithin{equation}{section}

\newcommand{\E}{\mathbf{E}}

\begin{document}

\title[Gaussian approximation of moments]{Gaussian approximation of moments of sums of independent random variables}

\author[M. Lis]{Marcin Lis}
\address{
Institute of Mathematics, University of Warsaw \\
Banacha 2, 02-097 Warszawa}\email{ml236089@students.mimuw.edu.pl}
\address{
Department of Mathematics, Vrije Universiteit Amsterdam \\
De Boelelaan 1081-1083, 1081HV Amsterdam}
\email{m.lis@few.vu.nl}

\begin{abstract}
We continue the research of Lata\l a \cite{latala} on improving estimates of $p$-th moments of sums of independent 
random variables. We generalize some of his results in the case when $2 \leq p \leq 4$ and present a combinatorial approach for even moments.
\end{abstract}



\maketitle

\section{Introduction}

In \cite{latala} Lata\l a studied moments of sums of independent random variables with logarithmically concave tails (i.e. variables for which the function $t\mapsto \ln(|X_i| \geq t)$ is concave from  $[0,\infty)$ to $[-\infty,0]$). Among other results he obtained approximations for moments which are described in Corollary \ref{corollary:Latala_1} and Corollary \ref{corollary:Latala_2}.
From these estimates one can conclude that moments of sums of independent variables with logarithmically concave tails are close to the corresponding moments of gaussian variables as long as the variances of separate variables are uniformly small.

In this paper we generalize some of these results to the case when the variables are no longer required to have logarithmically concave tails. This is achieved by means of Lemma \ref{lemma:ineq_char_func} which improves Lemma 1 from \cite{latala}, and by introducing a combinatorial approach which is described in Section~3. Following Lata\l a the first tool is used to obtain bounds for moments when $2 \leq p \leq 4$ and the combinatorial method gives estimates for moments of even degree. It also gives upper bounds for moments when the assumption of symmetry of the variables is dropped.

Throughout the paper we write $\|X\|_p=(\E |X|^p)^{1/p}$ for the $p$-th moment of the random variable $X$, and $\varphi_X$ for the characteristic function of $X$. By $\gamma_p$ we denote the $p$-th moment of a standard Guassian variable and by $(\varepsilon_k)_{1\leq k \leq n}$ a sequence of independent Rademacher variables, i.e. $\mathbf{P}(\varepsilon_k= \pm 1) =1/2$.

\section{Gaussian approximation for moments of order $2\leq p\leq4$}
We begin with the following simple lemma about characteristic functions. 
\begin{lem} 
If $X$ is a symmetric random variable, then
\begin{align*}
     1 - \frac{t^2}{2} \E X^2\leq \varphi_X (t) \leq 1 - \frac{t^2}{2}\E X^2 + \frac{t^4}{4!}\E X^4 . 
\end{align*}
\begin{proof} It follows immediately from the fact that $\varphi_X(t)= \E \cos tX$, and from the elementary inequality
\begin{align*}
1-\frac{x^2}{2}& \leq \cos x \leq 1 - \frac{x^2}{2}+ \frac{x^4}{4!} \text{ for all } x\in \mathbb{R}. \qedhere
\end{align*}
\end{proof}
\label{lemma:simple_char_func}
\end{lem}

The following lemma is the main ingredient of the results in the paper.
\begin{lem} 
\label{lemma:ineq_char_func}
Let $X_1, X_2, \ldots, X_n$ and $Y_1, Y_2, \ldots, Y_n $ be two sequences of independent symmetric random variables such that $\E X_k^2=\E Y_k^2$ for all $1 \leq k \leq n$. If $1\leq m < n$ is such that $\E X_l^2 = \max_{1\leq k \leq n} \E X_k^2$ for some $1\leq l \leq m$, and $\sum_{k=1}^m \E X_k^2 \geq \frac{1}{6} \max_{m < k \leq n} \frac{\E Y_k^4}{\E Y_k^2}$, then for $S=\sum_{k=1}^n X_k$ and $R=\sum_{k=m+1}^n Y_k$, it holds
\begin{align*}
\varphi_S(t) + \frac{t^2}{2}\sum_{k=1}^m \E X_k^2 \geq \varphi_R(t) \text{ for all } t \in \mathbb{R}. 
\end{align*}
\begin{proof} First consider the case when $r= \frac{t^2}{2} \sum_{k=1}^m \E X_k^2 < 1$. Then $\frac{t^2}{2}\E X_k^2 \leq \frac{t^2}{2}\E X_l^2 \leq r < 1$. Hence, by Lemma \ref{lemma:simple_char_func}, $\varphi_{X_k}(t), \varphi_{Y_k}(t) > 0$, and
\begin{align*}
\varphi_S(t)= \prod_{k=1}^n \varphi_{X_k}(t)\geq \prod_{k=1}^n(1-\frac{t^2}{2}\E X_k^2), \text{ and }
\end{align*}
\begin{align*}
\varphi_R(t)=\prod_{k=m+1}^n \varphi_{Y_k}(t) \leq \prod_{k=m+1}^n(1-\frac{t^2}{2}\E Y_k^2+\frac{t^4}{4!}\E Y_k^4).
\end{align*}
Thus, if we write $a_k=\frac{t^2}{2}\E X_k^2$, $r_k= \frac{t^2}{12}\frac{\E Y_k^4}{\E Y_k^2}$, it is enough to prove that 
\begin{align*}
\prod_{k=1}^n (1-a_n)+ \sum_{k=1}^m a_k \geq \prod_{k=m+1}^n(1-a_k+r_ka_k).
\end{align*}
Since $\prod_{k=1}^n (1-a_n) \geq (1-\sum_{k=1}^m a_k)\prod_{k=m+1}^n(1-a_k)$ and $0\leq r_k \leq r < 1 $, it suffices to prove that
\begin{align*}
(1-r) \prod_{k=m+1}^n(1-a_k) + r  \geq \prod_{k=m+1}^n(1-a_k+ra_k) \text{ for } 0\leq r \leq 1. 
\end{align*}
This is true since the left-hand side of the inequality is linear and the right-hand side is a convex function of $r$ on $[0,1]$, and the inequality holds true
for $r=0,1$.

If $r\geq1$ and $\varphi_S(t)> 0$, then the inequality holds true for obvious reasons.

If $\varphi_S(t)< 0$, then $\varphi_{X_i}(t) < 0$ for some $1\leq i \leq n$, and Lemma \ref{lemma:simple_char_func} implies that $|\varphi_{X_i}(t)|\leq \frac{t^2}{2}\E X_i^2-1$. Hence,
\begin{align*}
&\varphi_S(t)+\frac{t^2}{2}\sum_{k=1}^m \E X_k^2 \geq \frac{t^2}{2}\E X_l^2 - |\varphi_{X_i}(t)| \\
& \geq \frac{t^2}{2}(\E X_l^2-\E X_i^2)+1\geq 1 \geq \varphi_R(t).
\end{align*}
This finishes the proof.
\end{proof}
\end{lem}

\begin{rem} If we additionally assume that the sequence $\E X_k^2$, $k=1,\ldots,n$ is nonincreasing then 
$m=\big\lceil \frac{1}{6} \max_{1\leq k \leq n} \frac{\E Y_k^4}{(\E Y_k^2)^2} \big \rceil$ fulfills the assumptions of Lemma \ref{lemma:ineq_char_func} if only $m<n$.
\label{remark:simpler_m}
\end{rem}

\begin{rem} If a symmetric random variable $Y$ has a logarithmically concave tail, then by Proposition 1 of \cite{latala}, $\frac{1}{6} \frac{\E Y^4}{\E Y^2} \leq 1$. Therefore if $\E Y_1^2 = \max_{1\leq k \leq n} \E Y_k^2$, and $Y_2,\ldots,Y_n$ have logarithmically concave tails, then $m=1$ fulfills the assumptions of Lemma \ref{lemma:ineq_char_func}
 \label{remark:logconcave_m}
\end{rem}

The next result is a genereralization of Lemma 2 from \cite{latala} and its proof follows that in \cite{latala}.
\begin{lem}
Let $X_1, X_2, \ldots, X_n$, $Y_1, Y_2, \ldots, Y_n $ and $m$ be as in Lemma~\ref{lemma:ineq_char_func}. Then,
\begin{equation*}
       \E \Big |\sum_{k=1}^{n} X_k \Big |^p \geq \E \Big |\sum_{k=m+1}^{n}  Y_k \Big |^p \ \ \ \ \text{for} \  2\leq p \leq 4.
\end{equation*}
\label{lemma:ineq_moments_2_p_4}
\begin{proof}
Let $S= \sum_{k=1}^{n}  X_k$ and $R= \sum_{k=m+1}^{n}  Y_k $. It is enough to prove the above inequality for $2<p<4$. By Lemma 4.2 of  \cite{haagerup}, we have for any random variable $X$ with a finite fourth moment,
\begin{align*}
\E|X|^p = C_p \int_0^{\infty}\Big(\varphi_X(t)-1+\frac{1}{2}t^2\E|X|^2\Big)t^{-p-1}dt,
\end{align*}
where ${C_p=-\frac{2}{\pi}\sin(\frac{p\pi}{2})\Gamma (p+1) > 0}$.
By Lemma \ref{lemma:ineq_char_func} we have 
\begin{align*}
 \varphi_{S}(t) - \varphi_{R}(t) = \prod_{k=1}^n \varphi_{X_k}(t) - \prod_{k=m+1}^n \varphi_{Y_k} (t) \geq  - \frac{t^2}{2} \sum_{k=1}^{m} \E X_k^2,
\end{align*}
and thus 
\begin{align*}
 \E|S|^p - \E|R|^p &= C_p\int_0^{\infty}\Big(\varphi_{S}(t) - \varphi_{R}(t) + \frac{t^2}{2} \sum_{k=1}^{m} \E X_k^2\Big)t^{-p-1}dt \geq 0. \qedhere
\end{align*}
\end{proof}
\end{lem}

The following proposition and corollary give some estimates for moments of sums of independent random variables.
\begin{prop}
Let $X_1,X_2,\ldots,X_n$ be a sequence of independent symmetric random variables such that  $\E X_1^2= \max_{1\leq k \leq n} \E X_k^2$ and 
\begin{align}
m= \max_{1\leq k\leq n} \Big \lceil\frac{1}{6} \frac{\E X_k^4}{(\E X_k^2)^2} \Big \rceil < n.
\label{eqn:new_m}
\end{align}
Then, for $2 \leq p \leq 4$,
\begin{align*}
\gamma_p \Big(\sum_{k=2}^n \E X_k^2 \Big)^{1/2} \leq \Big \| \sum_{k=1}^n X_k \Big \|_p \leq \gamma_p\Big(\sum_{k=1}^n \E X_k^2 \Big)^{1/2} + \sqrt{3m} \|X_1\|_2.
\end{align*}

\label{proposition:moments_2_p_4}
\begin{proof} We may assume that the sequence $(\E X_k^2)_{1\leq k \leq n}$ is nonincreasing.
Let $(Y_k)_{1 \leq k \leq n}$ be a sequence of independent symmetric Gaussian random variables with $\E Y_k^2 = \E X_k^2$ for all $1\leq k \leq n$.
By Remark \ref{remark:logconcave_m}, $m=1$ fulfills the assumptions of Lemma \ref{lemma:ineq_moments_2_p_4} applied to the sequences $(X_k)_{1\leq k \leq n}$ and $(Y_k)_{1\leq k \leq n}$. Hence the lower bound follows. 
To show the upper bound we again use Lemma \ref{lemma:ineq_moments_2_p_4}, but with the roles of the sequences $(X_k)_{1\leq k \leq n}$, $(Y_k)_{1\leq k \leq n}$ interchanged and with $m$ as in \eqref{eqn:new_m} (see Remark \ref{remark:simpler_m}). 
\begin{align*}
 \Big \| \sum_{k=1}^n X_k \Big\|_p &- \gamma_p\Big(\sum_{k=1}^n \E X_k^2 \Big)^{1/2}  \leq \Big\|  \sum_{k=1}^n X_k \Big\|_p - \Big \| \sum_{k=m+1}^{n} X_k  \Big \|_p \\& \leq \Big \| \sum_{k=1}^{m} X_k \Big \|_p  \leq \Big \| \sum_{k=1}^{m} X_k \Big \|_4  \leq \sqrt{3m} \|X_1\|_2.
\end{align*}
The last inequality follows from the identity $\E\big(\sum_{k=1}^{m} X_k\big)^4=\sum_{k=1}^{m} \E X_k^4 + 6\sum_{k<l} \E X_k^2 \E X_l^2$.
\end{proof}
\end{prop}

\begin{cor}
Let $X_1, X_2, \ldots, X_n$ and $m$ be as in Proposition \ref{proposition:moments_2_p_4}. Then,
\begin{align*}
\Big | \Big \| \sum_{k=1}^n X_k \Big \|_p - \gamma_p\Big(\sum_{k=1}^n \E X_k^2 \Big)^{1/2} \Big | \leq \ \sqrt{3m} \| X_1 \|_{2} \ \ \text{for} \ 2 \leq p \leq 4.
\end{align*}
\label{corollary:gaussian_approx_2_p_4}
\begin{proof}
The upper bound follows from Proposition \ref{proposition:moments_2_p_4}. To prove the lower bound we use the lower estimate from Proposition \ref{proposition:moments_2_p_4}.
\begin{align*}
\Big \| \sum_{k=1}^nX_k \Big \|_p  & \geq \gamma_p \Big(\sum_{k=2}^n \E X_k^2 \Big)^{1/2} \geq \gamma_p\Big(\Big(\sum_{k=1}^n \E X_k^2 \Big)^{1/2} - (\E X_1^2)^{1/2}\Big)  \\& \geq \gamma_p\Big(\sum_{k=1}^n \E X_k^2 \Big)^{1/2} - 3^{1/4} (\E X_1^2)^{1/2}.  \qedhere
\end{align*}
\end{proof}
\end{cor}

These results allow us to improve some results of \cite{latala}. Corollary~3 from~\cite{latala} says that if $(X_k)_{1 \leq k \leq n}$ is a sequence of independent symmetric random variables 
with logarithmically concave tails, then for $p\geq 3$,
\begin{equation*}
\Big | \Big \| \sum_{k=1}^n X_k \Big \|_p - \gamma_p\Big(\sum_{k=1}^n \E X_k^2 \Big)^{1/2} \Big | \leq p \max_{1\leq k \leq n}\| X_k \|_2 .
\end{equation*}
\begin{cor}
\label{corollary:Latala_1}
The above statement 
holds true for $p\geq 2$.
\begin{proof}
It is enough to use Corollary \ref{corollary:gaussian_approx_2_p_4} and Remark \ref{remark:logconcave_m}. 
\end{proof}
\end{cor}
The same applies to Theorem 2 from \cite{latala} which states that if $(X_k)_{1\leq k \leq n}$ is like above, and moreover the sequence $(\E X_k^2)_{1\leq k \leq n}$ is nonincreasing, then for $p\geq3$,
\begin{align*}
\max \bigg \{ \gamma_p \Big( \sum_{k\geq \lceil p/2 \rceil} \E X_k^2 \Big)^{1/2}, \Big \| \sum_{k<p} X_k \Big \|_p \bigg \} & \leq \Big \| \sum_{k=1}^n  X_k \Big \|_p \\
& \leq \gamma_p \Big ( \sum_{k\geq \lceil p/2 \rceil} \E X_k^2 \Big)^{1/2} + \Big \| \sum_{k<p} X_k \Big \|_p.
\end{align*}
\begin{cor}
\label{corollary:Latala_2}
The above statement 
holds true for $p\geq2$.
\begin{proof}
For $p=2$ it is obvious. Let $2 < p\leq 3$. We have $\lceil p/2 \rceil = 2$ and the lower bound is a consequence of Proposition \ref{proposition:moments_2_p_4}.
The upper estimate can be obtained as follows
\begin{equation*}
\Big \| \sum_{k=1}^n  X_k \Big \|_p \leq   \Big \| \sum_{k \leq 2} X_k \Big \|_p + \Big \| \sum_{k>2} X_k \Big \|_p \leq  \Big \| \sum_{k < p}  X_k \Big \|_p + \gamma_p \Big( \sum_{k\geq \lceil p/2 \rceil} \E X_k^2 \Big)^{1/2},
\end{equation*}
where the last inequality follows from Lemma \ref{lemma:ineq_moments_2_p_4} applied to $(X_k)_{1\leq k \leq n}$ and a sequence of independent symmetric Gaussian variables having the same second moments as $(X_k)_{1\leq k \leq n}$, and Remark \ref{remark:logconcave_m}.
\end{proof}
\end{cor}

\section{Even moments of sums of symmetric random variables}
We will now use combinatorial methods to give bounds for even moments of sums of independent random variables.
In this section we will use only one property of symmetric random variables, namely if $X$ is symmetric and $\E|X|^{2r+1} < \infty$ for $r \in \mathbb{N}$ then $\E X^{2r+1} =0$. 

We will use the multi-index notation to simplify formulae appearing in our statements. An $n$-dimensional multi-index $\alpha$ is an $n$-tuple $(\alpha_1,\alpha_2,\ldots,\alpha_n)  \in \mathbb{N}^n_0$, where $\mathbb{N}_0$ is the set of non-negative integers. 
From now on, all appearing multi-indices will be $n$-dimensional. For a multi-index $\alpha$, we write $|\alpha|_i = \sum_{k=1}^i \alpha_k$,  $|\alpha|=|\alpha|_n$, $\alpha!=\prod_{k=1}^n\alpha_k!$, and $s(\alpha) = \{ k : \alpha_k \neq 0 \}$.
For a set $A$, we denote by $|A|$ its cardinality.

\begin{lem}
Let $r\in \mathbb{N}$, $r\geq 2$ and let $X_1, X_2, \ldots,X_n$ be a sequence of independent symmetric random variables such that $(\E X_k^2)_{1\leq k \leq n}$ is nonincreasing. Moreover, let $C \geq 1$ be such that
\begin{equation*}\E X_k^{2l} \leq C^{2l-2} \frac{(2l)!}{2^l} (\E X_k^2)^l \text{ for all } l,k\in \mathbb{N},\ 1 \leq l \leq r, \ 1\leq k \leq n,
\end{equation*} and $\lceil C^2(r-1) \rceil < n$. Then,
\begin{align*}
    \E \Big (\sum_{k=\lceil C^2(r-1) \rceil +1}^{n}  X_k \Big)^{2r} \leq \frac{(2r)!}{2^{r}}\mathop{\sum_{|\alpha| = r}}_{|s(\alpha)|=r} \prod_{k=1}^n (\E X_k^2)^{\alpha_k}.
\end{align*}
\begin{proof}
Let $D = \lceil C^2(r-1) \rceil+1$. We have
\begin{align*}
        \E \Big(\sum_{k=D}^{n} & X_k \Big)^{2r}  = \mathop{\sum_{|\alpha|=r}}_{\ |\alpha|_{D-1} = 0} \frac{(2r)!}{(2\alpha)!} \E X_{D}^{2\alpha_{D}} \E X_{D+1}^{2\alpha_{D+1}} \cdots \E X_n^{2\alpha_n} \\
        &\leq \mathop{\sum_{|\alpha|=r}}_{\ |\alpha|_{D-1} = 0} \frac{(2r)!}{(2\alpha)!}\prod_{k=1}^n (\E X_k^2)^{\alpha_k}
        \frac{(2\alpha_{D})! (2\alpha_{D+1})! \cdots (2\alpha_n)! }{2^{\alpha_{D}}  2^{\alpha_{D+1}} \cdots 2^{\alpha_n}}C^{2(r-|s(\alpha)|)} \\ 
        & = \frac{(2r)!}{2^{r}} \mathop{\sum_{|\alpha|=r}}_{\ |\alpha|_{D-1} = 0}C^{2(r-|s(\alpha)|)} \prod_{k=1}^n (\E X_k^2)^{\alpha_k}.
\end{align*} 
Thus it is enough to show that
\begin{equation*}
	\mathop{\sum_{|\alpha|=r}}_{\ |\alpha|_{D-1} = 0} C^{2(r-|s(\alpha)|)} \prod_{k=1}^n (\E X_k^2)^{\alpha_k} \leq \mathop{\sum_{|\alpha| = r}}_{|s(\alpha)|=r} \prod_{k=1}^n (\E X_k^2)^{\alpha_k},
\end{equation*}

If $I$ is a nonempty subset of $\{D,D+1,\ldots, n\}$ and $i=|I| \leq r$, then there are exactly $ \binom {r-1}{r-i}$ multi-indices $\alpha$ on the left-hand side which satisfy $s(\alpha) = I$ (because this is the number of different ways one can put $r$ indistinguishable balls in $i$ distinguishable urns, without leaving any urn empty). Moreover, there are exactly $\binom{D-1}{r-i}$ multi-indices $\alpha$ on the right-hand side such that $s(\alpha) \cap \{D,D+1,\ldots, n\} = I$, since we can add to these $i$ fixed indices any $r-i$ indices from the set $\{1,2,\ldots,D-1\}$. Since $(\E X_k^2)_{1\leq k \leq n}$ is nonincreasing, any term corresponding to the selected multi-indices on the right-hand side is larger than any chose term on the left-hand side. Thus it is enough to show that $\binom{D-1}{r-i} \geq C^{2(r-i)} \binom {r-1}{r-i}$. Indeed, we have 
\begin{align*}
\binom{D-1}{r-i}& /\binom {r-1}{r-i} = \frac{(D-1)(D-2)\cdots(D-r+i)}{(r-1)(r-2)\cdots(r-(r-i))} \\
&= \frac{\lceil C^2(r-1) \rceil} {r-1} \frac{\lceil C^2(r-1) \rceil-1}{r-2}\cdots \frac {\lceil C^2(r-1) \rceil-(r-i-1)}{i} \\
&\geq C^{2(r-i)}.
\end{align*}
We finish the proof by repeating this procedure for every $I \subseteq \{D,D+1,\ldots,n\}$.
\end{proof}
\label{lemma:comb_symmetric_moments}
\end{lem}
As an immediate consequence we get the following result.
\begin{lem} Let $X_1, X_2, \ldots,X_n$ and $C$ be as in Lemma \ref{lemma:comb_symmetric_moments}. Then,
\begin{equation*}
 \E \Big(\sum_{k=\lceil C^2(r-1) \rceil +1}^{n} X_k \Big)^{2r}  \leq \E \Big(\sum_{k=1}^{n} (\E X_k^2)^{1/2} \varepsilon_k \Big)^{2r}.
\end{equation*} 
\begin{proof}
We have
\begin{align*}
\E &\Big(\sum_{k=1}^{n}  (\E X_k^2)^{1/2} \varepsilon_k \Big)^{2r} 
  = \sum_{|\alpha|=r} \frac{(2r)!}{(2\alpha)!}\prod_{k=1}^n (\E X_k^2)^{\alpha_k} \\ &\geq 	\frac{(2r)!}{2^{r}}\mathop{\sum_{|\alpha| = r}}_{|s(\alpha)|=r} \prod_{k=1}^n (\E X_k^2)^{\alpha_k} \geq 	\E \Big(\sum_{k=\lceil C^2(r-1) \rceil +1}^{n}  X_k \Big)^{2r}.
\end{align*}
where the last inequality follows from Lemma \ref{lemma:comb_symmetric_moments}.
\end{proof}
\label{corollary:symmetric_even_moments}
\end{lem} 
If $\mathcal{E}$ is a random variable with symmetric exponential distribution with variance $1$, then $\E\mathcal{E}^{2l} = (2l)!/2^l$. Therefore
if $X_1, X_2, \ldots,X_n$ have logarithmically concave tails then $C=1$ by Proposition~1 from \cite{latala}, and we get a different proof 
of Theorem 1 from \cite{latala} for $p=2r$.
From the corollary above we can conclude the following results.
\begin{prop} 
Let $X_1, X_2, \ldots,X_n$ and $C$ be as in Lemma \ref{lemma:comb_symmetric_moments}. Then,
\begin{equation*}
   			\gamma_{2r}\Big(\sum_{k=r}^n \E X_k^2\Big)^{1/2} \leq \Big\| \sum_{k=1}^n X_k \Big\|_{2r} \leq \gamma_{2r}\Big(\sum_{k=1}^n \E X_k^2\Big)^{1/2}  + 2\lceil C^2(r-1) \rceil \|X_1\|_{2}.
\end{equation*} 
\begin{proof}
The lower bound follows from Lemma \ref{corollary:symmetric_even_moments} applied to a sequence of independent symmetric Gaussian variables with variances $(\E X_k^2)_{1\leq k  \leq n}$ (where $C=1$), and from the fact that $\big\| \sum_{k=1}^n (\E X_k^2)^{1/2}\varepsilon_k \big\|_{2r} \leq \big\| \sum_{k=1}^n X_k \big\|_{2r} $. Let $D=\lceil C^2(r-1) \rceil$. We have 
\begin{align*}
\Big\| \sum_{k=1}^n X_k \Big \|_{2r} &- \gamma_{2r}\Big(\sum_{k=1}^n \E X_k^2\Big)^{1/2} \leq \Big\| \sum_{k=1}^n X_k \Big\|_{2r} - \Big\| \sum_{k=D+1}^n X_k \Big\|_{2r} \leq \Big\| \sum_{k=1}^{D} X_k \Big\|_{2r} 
\\ & \leq  \Big\| \sum_{k=1}^{2D} \sigma_k^{\circ} \varepsilon_k \Big\|_{2r} \leq
2D \|X_1\|_{2 },
\end{align*}
where $\sigma_k^{\circ} = (\E X_1^2)^{1/2}$ for $1\leq k \leq D$, and $\sigma_k^{\circ} = (\E X_{k-D}^2)^{1/2}$ for $D < k \leq 2D$.
We used here the fact that $\gamma_{2r}\big(\sum_{k=1}^n \E X_k^2\big)^{1/2} \geq \big \| \sum_{k=1}^n (\E X_k^2)^{1/2}
\varepsilon_k \big \|_{2r}$, and Lemma \ref{corollary:symmetric_even_moments} twice. Hence, the upper bound holds true.
\end{proof}
\label{corollary:symmetric_even_estimation}
\end{prop}

\begin{cor} 
Let $X_1, X_2, \ldots,X_n$ and $C$ be as in Lemma \ref{lemma:comb_symmetric_moments}. Then, 
\begin{equation*}
   			\Big | \Big\| \sum_{k=1}^n X_k \Big\|_{2r} - \gamma_{2r}\Big(\sum_{k=1}^n \E X_k^2\Big)^{1/2} \Big | \leq 2\lceil C^2(r-1) \rceil \|X_1\|_{2}.
\end{equation*} 
\begin{proof}
The upper bound follows from Proposition \ref{corollary:symmetric_even_estimation}. To prove the lower bound we use the lower estimate from Proposition \ref{corollary:symmetric_even_estimation}, and Lemma~\ref{corollary:symmetric_even_moments},
\begin{align*}
\Big\| \sum_{k=1}^n X_k \Big\|_{2r} & \geq 	\gamma_{2r}\Big(\sum_{k=r}^n \E X_k^2 \Big)^{1/2} \geq \gamma_{2r}\Big(\sum_{k=1}^n \E X_k^2\Big)^{1/2} -\gamma_{2r}\Big(\sum_{k=1}^{r-1} \E X_k^2 \Big)^{1/2} \\
&\geq \gamma_{2r}\Big(\sum_{k=1}^n \E X_k^2\Big)^{1/2} -\Big\|\sum_{k=1}^{2(r-1)} \sigma_k^{\circ}\varepsilon_k\Big\|_{2r}\\
 & \geq \gamma_{2r}\Big(\sum_{k=1}^n \E X_k^2\Big)^{1/2} - 2(r-1)\|X_1\|_{2},
\end{align*}
where $\sigma_k^{\circ} = (\E X_1^2)^{1/2}$ for $1\leq k \leq r-1$, and $\sigma_k^{\circ} = (\E X_{k-r+1}^2)^{1/2}$ for $r-1 < k \leq 2r-2$.
\end{proof}
\label{corollary:symmetric_even_gaussian_approx}
\end{cor}

\section{Even moments of sums of centered random variables}
We begin with a lemma which is very similar to Lemma \ref{lemma:comb_symmetric_moments} but since our variables are no longer symmetric we get a weaker conclusion.
\begin{lem} 
\label{lemma:comb_nonsymmetric_moments}
Let $r\in \mathbb{N}$, $r\geq 2$ and let $X_1, X_2, \ldots,X_n$ be a seqence of independent centered random variables such that $(\E X_k^2)_{1\leq k \leq n}$ is nonincreasing. Moreover, let $C \geq 1$ be such that 
\begin{equation*}
|\E X_k^{l}| \leq C^{l-2}\frac{l!}{2^{l/2}} (\E X_k^2)^{l/2} \text{ for all } l,k\in \mathbb{N},\ 2 \leq l \leq 2r,\ 1\leq k \leq n 
\end{equation*} and $\lceil  C^2\frac{r(r-1)}{2} \rceil < n$. Then,
\begin{align*}
    \E \Big (\sum_{k=\lceil C^2\frac{r(r-1)}{2} \rceil +1}^{n}  X_k \Big)^{2r} \leq \frac{(2r)!}{2^{r}}\mathop{\sum_{|\alpha| = r}}_{|s(\alpha)|=r} \prod_{k=1}^n (\E X_k^2)^{\alpha_k}. 
\end{align*}
\begin{proof}
Let $D = \lceil C^2\frac{r(r-1)}{2} \rceil+1$. We have
\begin{align*}
        \E \Big(\sum_{k=D}^{n} & X_k \Big)^{2r}  = \mathop{\sum_{|\alpha|=2r,\ \text{sing}(\alpha)=\emptyset}}_{\ |\alpha|_{D-1} = 0} \frac{(2r)!}{\alpha!}  \E X_{D}^{\alpha_{D}} \E X_{D+1}^{\alpha_{D+1}} \cdots \E X_n^{\alpha_n} \\
        &\leq \mathop{\sum_{|\alpha|=2r,\ \text{sing}(\alpha)=\emptyset}}_{\ |\alpha|_{D-1} = 0} \frac{(2r)!}{\alpha!}\prod_{k=1}^n (\E X_k^2)^{\alpha_k/2} 
        \frac{\alpha_{D}! \alpha_{D+1}! \cdots \alpha_n! }{2^{\alpha_{D}/2}  2^{\alpha_{D+1}/2} \cdots 2^{\alpha_n/2}}C^{2(r-|s(\alpha)|)} \\ 
        & = \frac{(2r)!}{2^{r}} \mathop{\sum_{|\alpha|=2r,\ \text{sing}(\alpha)=\emptyset}}_{\ |\alpha|_{D-1} = 0}C^{2(r-|s(\alpha)|)} \prod_{k=1}^n (\E X_k^2)^{\alpha_k/2},
\end{align*} 
where $\text{sing} (\alpha) = \{ k: \alpha_k=1 \}$.
Thus it is enough to show that
\begin{equation*}
	\mathop{\sum_{|\alpha|=2r,\ \text{sing}(\alpha)=\emptyset}}_{\ |\alpha|_{D-1} = 0} C^{2(r-|s(\alpha)|)} \prod_{k=1}^n (\E X_k^2)^{\alpha_k/2} \leq \mathop{\sum_{|\alpha| = r}}_{|s(\alpha)|=r}\prod_{k=1}^n (\E X_k^2)^{\alpha_k}.
\end{equation*}
Again, if $I$ is a nonempty subset of $\{D,D+1,\ldots, n\}$ and $|I| = i \leq r$, then there are exactly $\binom {2r-i-1}{2(r-i)}$ multi-indices $\alpha$ on left-hand side which satisfy $s(\alpha) = I$ (because this is the number of different ways one can put $2r$ indistinguishable balls in $i$ distinguishable urns in such a way that in each urn there are at least two balls). Moreover, the number of multi-indices $\alpha$ for which $s(\alpha) \cap \{D,D+1,\ldots, n\} = I$ is $\binom{D-1}{r-i}$. We notice again that any term corresponding to the selected multi-indices on the right-hand side is larger than any chosen on the left-hand side. Thus it is enough to show that 
$\binom{D-1}{r-i} \geq C^{2(r-i)} \binom {2r-i-1}{2(r-i)}$. Indeed, we have 
\begin{align*}&\binom{D-1}{r-i}/\binom {2r-i-1}{2(r-i)} = \frac{(2(r-i))!}{(r-i)!}\frac{(D-1)(D-2)\cdots(D-(r-i))}{(2r-i-1)\cdots(i+1)i} \\
&= \frac{(2(r-i))!}{(r-i)!}  \frac{\lceil C^2r(r-1)/2 \rceil-(r-i-1)} {i(2r-i-1)} \\ & \ \ \ \ \ \ \ \ \ \ \ \ \ \ \ \ \cdot\frac{\lceil C^2r(r-1)/2 \rceil-(r-i-2)}{(i+1)(2r-i-2)}\cdots \frac {\lceil C^2r(r-1)/2 \rceil}{(r-1)r} \\
&\geq C^{2(r-i)}\frac{r(r-1)-2(r-i-1)} {i(2r-i-1)} \frac{ r(r-1)-2(r-i-2)}{(i+1)(2r-i-2)}\cdots \frac {r(r-1)}{(r-1)r}.
\end{align*}
 Therefore it suffices to prove that $r(r-1)-2(r-j-1) \geq j(2r-j-1)$ for $1\leq j\leq r-1$. By substituting $j$ for $r-j$ we need to prove that $j^2-3j+2\geq 0$, which is true for $j \in \mathbb{N}_+$.
We finish the proof by repeating this procedure for every $I \subseteq \{D,D+1,\ldots,n\}$. 
\end{proof}
\end{lem}

We get as a consequence the following results. 
\begin{lem} Let $X_1, X_2, \ldots,X_n$ and $C$ be as in Lemma \ref{lemma:comb_nonsymmetric_moments}. Then,
\begin{equation*}
   				\E\Big(\sum_{k=\lceil C^2\frac{r(r-1)}{2} \rceil+1}^{n}  X_k \Big)^{2r} \leq \E \Big(\sum_{k=1}^{n} (\E X_k^2)^{1/2} \varepsilon_k \Big)^{2r}.
\end{equation*} 
\label{corollary:nonsymmetric_even_moments}
\end{lem}

\begin{prop} Let $X_1, X_2, \ldots,X_n$ and $C$ be as in Lemma \ref{lemma:comb_nonsymmetric_moments}. Then,
\begin{equation*}
   		 \Big\| \sum_{k=1}^n X_k \Big\|_{2r} \leq \gamma_{2r}\Big(\sum_{k=1}^n \E X_k^2\Big)^{1/2} + 2 \big \lceil C^2\frac{r(r-1)}{2} \big \rceil \|X_1\|_{2}.
\end{equation*} 

\begin{proof}
Let $D= \lceil C^2\frac{k(k-1)}{2}  \rceil$. We have 
\begin{align*}
\Big\|\sum_{k=1}^n X_k \Big\|_{2r} &- \gamma_{2r}\Big(\sum_{k=1}^n \E X_k^2\Big)^{1/2} \leq \Big\|\sum_{k=1}^n X_k \Big\|_{2r} - \Big\| \sum_{k=D+1}^n X_k \Big\|_{2r} \\
&  \leq \Big\| \sum_{k=1}^{D} X_k \Big\|_{2r} \leq \Big\| \sum_{k=1}^{2D} \sigma_k^{\circ}\varepsilon_k \Big\|_{2r}  \leq 2D \|X_1\|_{2},
\end{align*}
where $\sigma_k^{\circ} = (\E X_1^2)^{1/2}$ for $k \leq D$, and $\sigma_k^{\circ} = (\E X_{k-D}^2)^{1/2}$ for $D < k \leq 2D$. We used here Lemma \ref{corollary:nonsymmetric_even_moments} twice and again the fact that $\gamma_{2r}\Big(\sum_{k=1}^n \E X_k^2\Big)^{1/2} \geq \big\| \sum_{k=1}^n (\E X_k^2)^{1/2}
\varepsilon_k \big \|_{2r}$.
\end{proof}
\end{prop}

\section{Arbitrary $p$-th moments for $p\geq 4$}
 For a multi-index $\alpha$, 
we write $\sigma^{2\alpha} = \prod_{k=1}^n \sigma_k^{2\alpha_k}$.

\begin{lem} \label{big_lemma_1}
If $\sigma \in \mathbb{R}^n$ is such that $|\sigma_1| \geq |\sigma_2| \geq \ldots \geq |\sigma_n|$. Then for all $r \geq 1$,
\begin{align*}
   \frac{2r+1}{2r-1} \bigg(\E \Big(\sum_{k=1}^{n} \sigma_k \varepsilon_k \Big)^{2r}\bigg)^2 \geq \frac{(2r+2)!}{2^{r+1}}  \bigg( \mathop{\sum_{|\alpha| = r+1}}_{|s(\alpha)|=r+1}  \sigma^{2\alpha} \bigg) \E \Big(\sum_{k=1}^{n} \sigma_k \varepsilon_k \Big)^{2r-2}.
    \end{align*}

\begin{proof}
We have
\begin{align*}
        \bigg(\E \Big(\sum_{k=1}^{n} \sigma_k \varepsilon_k \Big)^{2r}\bigg)^2  &= \bigg(\sum_{|\alpha|=r}\frac{(2r)!}{(2\alpha)!}\sigma^{2\alpha}\bigg)^2 
          = (2r!)^2\sum_{|\alpha|=|\beta|=r}\frac{\sigma^{2(\alpha+\beta)}}{(2\alpha)!(2\beta)!},
\end{align*} 
and 
\begin{align*}
   \frac{(2r+2)!}{2^{r+1}}    \bigg( \mathop{\sum_{|\alpha| = r+1}}_{|s(\alpha)|=r+1}  & \sigma^{2\alpha} \bigg) \E \Big(\sum_{k=1}^{n} \sigma_k \varepsilon_k \Big)^{2r-2} \\ & = (2r+2)!(2r-2)! \mathop{\sum_{|\alpha| = r+1,|\beta|=r-1}}_{|s(\alpha)|=r+1}  \frac{\sigma^{2(\alpha+\beta)}}{2^{r+1}(2\beta)!}.  
\end{align*} 
We define two subsets of $\mathbb{N}^n_0\times \mathbb{N}^n_0$ as follows:
\begin {align*}
 L&=\{(\alpha,\beta): |\alpha|=|\beta|=r \}, \\ 
 R&=\{(\gamma,\delta):|\gamma|=r-1,\ |\delta| = r+1,\ |s(\delta)|=r+1\}.
\end{align*}
Since $ \frac{2r+1}{2r-1} (2r)!^2 = \frac{r}{r+1} (2r+2)!(2r-2)! $ to prove the desired inequality, we need to show that
\begin{align} \label{eqn:big_lemma_2}
 	\frac{r}{r+1} \sum_{(\alpha,\beta)\in L}\frac{\sigma^{2(\alpha+\beta)}}{(2\alpha)!(2\beta)!} \geq  \sum_{(\gamma,\delta)\in R} \frac{\sigma^{2(\gamma+\delta)}}{2^{r+1}(2\gamma)!}.
\end{align} 
To prove this, we will divide $R$ into disjoint subsets. For each such subset we will find a corresponding subset of $L$. We will make sure that the subsets of $L$ are also disjoint and that the sums over corresponding subsets satisfy the desired inequality. For $(\gamma,\delta)\in R$, we define 
\begin{align*}
R(\gamma,\delta) = \{ (\gamma',\delta') :\ & \gamma'+\delta' = \gamma+\delta, \\ &
															 \gamma'_k=\gamma_k \ \text{and} \ \delta'_k=\delta_k \ \text{for} \ k \notin \text{sing}(\gamma+\delta), \\
															& |\gamma'| = r-1, \ |\delta'| =r+1, \\ &
															 |s(\delta')|=r+1 \} \subseteq R
\end{align*}
and 
\begin{align*}
L(\gamma,\delta) = \{ (\alpha,\beta) :\ & \alpha+\beta = \gamma+\delta,  \\&
															 \alpha_k=\gamma_k \ \text{and} \ \beta_k=\delta_k \ \text{for} \ k \notin \text{sing}(\gamma+\delta), \\
															& |\alpha| =  |\beta| =r\} \subseteq L.
\end{align*}
One can easily see that both families $\mathcal{R}=\{R(\gamma,\delta):(\gamma,\delta)\in R\}$ and $\mathcal{L}=\{L(\gamma,\delta):(\gamma,\delta)\in R\}$ are pairwise disjoint. Moreover $\mathcal{R}$ is a partition of $R$. It is also not difficult to notice that 
$R(\gamma,\delta)=R(\gamma',\delta')$ if and only if $L(\gamma,\delta)=L(\gamma',\delta')$. From the definition of $R(\gamma,\delta)$ and $L(\gamma,\delta)$, it follows also that the function $(\alpha,\beta) \mapsto \frac{\sigma^{2(\alpha+\beta)}}{(2\alpha)!(2\beta)!}$ is constant on 
$R(\gamma,\delta)$ and $L(\gamma,\delta)$ and takes the same value on both sets. Thus to prove \eqref{eqn:big_lemma_2} it is enough to show that 
\begin{align*}
	\frac{r}{r+1} |L(\gamma,\delta)| \geq |R(\gamma,\delta)| \ \text{for every} \ (\gamma,\delta) \in R.
\end{align*} 
We say that $k$ is \emph{single} in $\alpha$ if $k \in \text{sing}(\alpha)$.
Let us fix $(\gamma,\delta) \in R$. We have
\begin{align*}
s(\gamma) \setminus (\text{sing}(\gamma) \cap \text{sing}(\gamma+\delta)) \supseteq s(\delta) \setminus (\text{sing}(\delta) \cap \text{sing}(\gamma+\delta)),
\end{align*}
since all indices in $s(\delta)$ are single in $\delta$ and these indices from $s(\delta)$ which are not single in $\gamma+\delta$ must appear in $s(\gamma)$. Thus we have
\begin{align*}
r-1- |\text{sing}(\gamma) \cap \text{sing}(\gamma+\delta)| & \geq |s(\gamma) \setminus (\text{sing}(\gamma) \cap \text{sing}(\gamma+\delta))|  \\
& \geq |s(\delta) \setminus (\text{sing}(\delta) \cap \text{sing}(\gamma+\delta))| \\
& = r+1 - |\text{sing}(\delta) \cap \text{sing}(\gamma+\delta)|,
\end{align*} 
and therefore,
\begin{align}
|\text{sing}(\delta) \cap \text{sing}(\gamma+\delta)| \geq |\text{sing}(\gamma) \cap \text{sing}(\gamma+\delta)|+2.
    \label{eqn:big_lemma_3}
\end{align}
To simplify notation we define $g = |\text{sing}(\gamma) \cap \text{sing}(\gamma+\delta)|$ and $d = |\text{sing}(\delta) \cap \text{sing}(\gamma+\delta)|$. The idea behind the definition of the set $R(\gamma,\delta)$ is that it contains elements resulting from $(\gamma,\delta)$ by replacing the original sets $\text{sing}(\gamma) \cap \text{sing}(\gamma+\delta)$ and $\text{sing}(\delta) \cap \text{sing}(\gamma+\delta)$ by new subsets of $\text{sing}(\gamma+\delta)$ while preserving the cardinality of these sets. Thus $ |R(\gamma,\delta)| = \binom{d+g}{d}$. We can interpret $L(\gamma,\delta)$ similarly but this time we change the cardinality of the sets 
of single indices in a proper way, and we get $ |L(\gamma,\delta)| = \binom{d+g}{d-1}$.
From \eqref{eqn:big_lemma_3} we get  
$
|L(\gamma,\delta)|=\binom{d+g}{d-1} = \frac{d}{g+1}\binom{d+g}{d} \geq \frac{g+2}{g+1}\binom{d+g}{d} \geq \frac{r+1}{r} \binom{d+g}{d} =\frac{r+1}{r} |R(\gamma,\delta)|,
$ because $g+2\leq d\leq r+1$,
which ends the proof.
\end{proof}
\label{lemma:big_lemma}
\end{lem}

\begin{rem}
The constant $\frac{2r+1}{2r-1}$ in the above lemma is optimal. It is enough to take all $\sigma_k = 1/\sqrt{n}$ 
and see that both sides of the inequality approach the same limit as $n$ goes to infinity.
\end{rem}

We can now prove an inequality for moments of arbitrary order which is weaker than previous results for even moments because of its non-unital multiplicative constant. 

\begin{prop}
Let $r\in \mathbb{N}$, $r\geq1$ and let $X_1, X_2, \ldots,X_n$ be a sequence of independent centered random variables such that $(\E X_k^2)_{1\leq k \leq n}$ is nonincreasing. Moreover, let $C \geq 1$ be such that 
\begin{equation*}|\E X_k^{l}| \leq C^{l-2}\frac{l!}{2^{l/2}} (\E X_k^2)^{l/2}\text{ for all } l,k\in \mathbb{N},\ 2 \leq l \leq 2r,\ 1\leq k \leq n.
\end{equation*}
Then, for $2\leq p\leq 2r$,
\begin{equation}
	 \frac{2 \lfloor p/2 \rfloor +1 }{2 \lfloor p/2 \rfloor -1} \E \Big|\sum_{k=1}^{n} (\E X_k^2)^{1/2} \varepsilon_k \Big|^p \geq  \E \Big|\sum_{k=C_p}^{n}  X_k \Big|^p
    \label{eqn:general}
\end{equation} if only $C_p \leq n$,
where $C_p=\lceil C^2 \lfloor p/2 \rfloor \rceil +1 $ when all $X_k$ are symmetric, and $C_p=\lceil C^2 \frac{\lfloor p/2 \rfloor (\lfloor p/2 \rfloor+1)}{2} \rceil +1 $ otherwise.
\begin{proof}
We will use the fact that for every random variable $X$, the function $p \mapsto \log \E |X|^p$ is convex. We define 
\begin{align*}
f(p) = \log \bigg [\E \Big |\sum_{k=1}^{n} (\E X_k^2)^{1/2} \varepsilon_k \Big |^p \bigg ] \ \ \text{and} \ \ 
 g(p) = \log \bigg [\frac{2 \lfloor p/2 \rfloor -1 }{2 \lfloor p/2 \rfloor +1} \E \Big |\sum_{k= C_p }^{n}  X_k \Big |^p \bigg]. 
\end{align*}
Note that $f$ is convex on $(0,\infty)$ and $g$ is convex on $(2l,2l+2)$ for every $1\leq l \leq r$, since $\lfloor p/2 \rfloor$ is constant on such intervals. Let us fix  $l$. By convexity of $g$ on $(2l,2l+2)$, we can identify $g$ with its continuous extension to the closed interval $[2l,2l+2]$. We will show that \eqref{eqn:general}  holds on $[2l,2l+2]$. Let $p\in[2l,2l+2]$. Since $f$ is convex, we have
\begin{equation*}
f \geq f(2l) + \frac{p-2l}{2}\big(f(2l) - f(2l-2)\big) =: \underline{f}(p),
\end{equation*}
and since $g$ is convex on $[2l,2l+2]$, we have
\begin{equation*}
g \leq g(2l) + \frac{p-2l}{2}\big(g(2l+2) - g(2l)\big) =: \overline{g}(p).
\end{equation*} 

Therefore it is enough to show that $\underline{f}(p)  \geq \overline{g}(p)$. By Lemma \ref{corollary:symmetric_even_moments} and~\ref{corollary:nonsymmetric_even_moments}, 
$\underline{f}(2l) = f(2l) \geq g(2l) = \overline{g}(2l)$, and because of linearity of both $\underline{f}$ and $\overline{g}$, we only need 
to show that $\underline{f}(2l+2) = 2 \log \E |\sum_{k=1}^{n} (\E X_k^2)^{1/2}\varepsilon_k |^{2l} - \log \E |\sum_{k=1}^{n} (\E X_k^2)^{1/2} \varepsilon_k |^{2l-2}  \geq \log \frac{2l-1}{2l+1} \E |\sum_{k=C_{2l}}^{n} X_k |^{2l+2} = \overline{g}(2l+2)$. This follows from Lemma \ref{corollary:symmetric_even_moments}, Lemma \ref{corollary:nonsymmetric_even_moments} and Lemma \ref{lemma:big_lemma}.
\end{proof}
\end{prop}

\section*{Aknowledgments}
The author would like to thank Professor Rafa\l\ Lata\l a for his guidance
through the subject and Stanis\l aw Kwapie\'n for his help with the editting of the paper.

\end{document}